
\input amstex
\loadbold


\font\tenbsy=cmbsy10
\font\sevenbsy=cmbsy7
\font\fivebsy=cmbsy5

\skewchar\tenbsy='60 
\skewchar\sevenbsy='60 
\skewchar\fivebsy='60

\newfam\bsyfam

\textfont\bsyfam=\tenbsy
\scriptfont\bsyfam=\sevenbsy
\scriptscriptfont\bsyfam=\fivebsy

\define\aihp{Ann.\ Inst.\ H.\ Poincar\'e Probab.\ Statist.}
\define\jsp{J.\ Statist.\ Phys.}

\define\ap{Ann.\ Probab.}

\define\ptrf{Probab.\ Theory Related Fields}

\define\spa{Stochastic Process.\ Appl.}

\define\superzd{^{\raise1pt\hbox{$\scriptstyle {\Bbb Z}^d$}}}
\define\xsubvn{x_{\lower1pt\hbox{$\scriptstyle V(n)$}}}
\define\xssubvn{x_{\lower1pt\hbox{$\scriptscriptstyle V(n)$}}}

\define\ssub#1{_{\lower1pt\hbox{$\scriptstyle #1$}}} 
\define\ssubb#1{_{\lower2pt\hbox{$\scriptstyle #1$}}} 
\define\Ssub#1{_{\lower1pt\hbox{$ #1$}}} 
\define\Ssu#1{_{\hbox{$ #1$}}} 
\define\Ssubb#1{_{\lower2pt\hbox{$ #1$}}} 
\define\sssub#1{_{\lower1pt\hbox{$\scriptscriptstyle #1$}}} 
\define\sssubb#1{_{\lower2pt\hbox{$\scriptscriptstyle #1$}}} 
\define\sssu#1{_{\hbox{$\scriptscriptstyle #1$}}}                                              
                                     
\define\<{\langle}
\define\>{\rangle}

\define\({\left(}
\define\){\right)}
\define\[{\left[}                                   
\define\]{\right]}

\define\lbrakk{\biggl\{}
\define\rbrakk{\biggr\}}


\define\osup#1{^{\raise1pt\hbox{$ #1 $}}}  
\define\osupm#1{^{\raise1pt\hbox{$\mskip 1mu #1 $}}}  
\define\osupp#1{^{\raise2pt\hbox{$ #1 $}}}  
\define\osuppp#1{^{\raise3pt\hbox{$ #1 $}}}  
\define\sosupp#1{^{\raise2pt\hbox{$\ssize #1 $}}}  
\define\sosup#1{^{\raise1pt\hbox{$\scriptstyle \mskip 1mu #1 $}}} 
\define\ssosup#1{^{\raise1pt\hbox{$\scriptscriptstyle #1 $}}} 
\define\ssosupp#1{^{\raise2pt\hbox{$\scriptscriptstyle #1 $}}} 
\define\ssosuppp#1{^{\raise3pt\hbox{$\scriptscriptstyle #1 $}}}



\define\nutil{{\tilde \nu}}
\define\etatil{{\tilde \eta}}


\define\e{\varepsilon}
\define\a{\alpha}                   
  
\define\ind{\text{I}}

\define\epi#1{\text{epi}\,#1}

\define\pikkuhyppy{\vskip .1in}


\define\mmR{{\bold R}}

\define\mmZ{{\bold Z}}

\define\mmal{{\boldsymbol \alpha}}


\define\SS{{\Cal S}}
\define\cAA{{\Cal A}}    

\define\CC{{\Cal C}}
\define\DD{{\Cal D}}

\define\FF{{\Cal F}}

\define\NN{{\Cal N}}

\define\TT{{\Cal T}}



\define\bP{{\Bbb P}}

\documentstyle{amsppt}
\magnification=\magstep1
\document     
\baselineskip=12pt
\pageheight{43pc} 

\NoBlackBoxes

\centerline{\bf Recent Results and Open Problems on the
Hydrodynamics}

\pikkuhyppy

\centerline{\bf  of Disordered Asymmetric Exclusion and Zero-Range
Processes}

\pikkuhyppy

$$\text{1998}$$

\hbox{}

\centerline{  Timo  Sepp\"al\"ainen
\footnote""{ Research partially supported by NSF grant DMS-9801085. }}
\hbox{}
\centerline{Department of Mathematics}
\centerline{Iowa State University}
\centerline{Ames, Iowa 50011, USA}
\centerline{seppalai\@iastate.edu}
\vfil

\flushpar
{\it Summary.} This paper summarizes 
results and some open problems  about the large-scale and long-time
behavior of  asymmetric, 
disordered exclusion and zero-range
processes. These processes have randomly chosen jump rates 
at the sites of the underlying lattice $\mmZ^d$. 
The interesting feature is that for 
suitably distributed random rates there is a phase
transition where the process behaves differently 
at high and low densities. Some of this distinction is 
visible on the hydrodynamic scale. But to fully understand 
the  phase
transition, results on a finer scale are needed. 

\vfil

\flushpar
Mathematics Subject Classification: Primary
 60K35, Secondary  82C22 

\flushpar
Keywords: Exclusion process, zero-range process,
 hydrodynamic limit, quenched disorder, phase transition

\flushpar
Short Title: Hydrodynamics for disordered particle systems 

\break

\head 1.  Introduction \endhead
This paper introduces some
recent results on the hydrodynamics of disordered, 
asymmetric simple exclusion processes
(SEP) and zero-range processes (ZRP). The
{\it disorder} refers to the rates of jumping attached to 
the sites of the underlying lattice: The particles 
move on $\mmZ^d$, and each site $x\in\mmZ^d$ has a 
random variable $\a_x$ that influences the exponential
rate at which particles leave site $x$. In SEP
 $\a_x$ is exactly the rate of jumping from $x$,
and in ZRP $\a_x$ multiplies the rate $r(\eta(x))$ that
depends on the number $\eta(x)$
 of particles currently occupying
site $x$. The {\it asymmetry} pertains to the jump probabilities
$p(x,y)$, according to which a particle jumping from $x$ 
chooses its new location $y$. We assume throughout that
the kernel $p(x,y)$ is translation invariant so that
$p(x,y)=p(0,y-x)\equiv p(y-x)$. Asymmetric jumping
means that there typically is a drift:
$\gamma\equiv\sum xp(x)\ne 0$. The assumption $\gamma\ne 0$
is not always  necessary, but without it
the limiting macroscopic conservation law 
becomes trivial. Some theorems require a stronger
assumption of 
{\it total asymmetry}: the dimension $d=1$, and all jumps
proceed to the right: $p(1)=1$. 

The disorder can also be attached to particles,
so that individual particles carry their own 
randomly chosen jump rates. We do not explicitly consider
such processes. 
One special case,  the totally asymmetric
simple exclusion process (TASEP) with particlewise
disorder, is partially covered by our 
discussion. This is  because  the gaps between the exclusion
particles with random rates can be regarded as the occupation
numbers of a ZRP
with random rates on the sites. This special
case has been studied in the physics literature as a 
model for traffic. See Krug (1998), Krug and Ferrari (1996), and 
their references. 

 The interesting 
phenomenon that appears in   disordered particle systems is 
a phase transition where  the process behaves differently 
at high and low densities. It occurs  when the 
distribution of the random rate
 has a sufficiently thin tail at its left endpoint $c>0$. 
Not much rigorous mathematical work exists
on this phase transition. 

A brief overview of the paper: In Section 2 the disordered
ZRP is described, together with two theorems.
 In Section
3 the same is done for the disordered SEP. Section 4
lists four open problems. Section 5 contains some proofs
and some comments on proofs. In particular, we included
in Section 5 a rigorous construction of a disordered
ZRP in $\mmZ^d$, and a proof of the invariance of 
a certain family of product measures. 
The construction is based on the percolation 
approach of Harris (1972). 

The hope is that this paper would be at least partially
accessible to the nonexpert. This is the motivation 
for inclusion of the proofs in Section 5, which are 
often referred to but less often spelled out in the literature.
For the same reason an attempt has been made to employ
 precise and complete  notation. 
This may make the text somewhat heavy to follow 
at times, but the alternative is to risk confusing the reader
who is not well-acquainted with disordered
particle systems. Of course, such an outcome may be
unavoidable in any case. 

Some familiarity with the 
 subject of interacting particle systems is
required for reading this paper. General references 
on particle systems 
are  Durrett (1988, 1995), 
Griffeath (1979), and Liggett (1985).
 References on hydrodynamical
limits are lectures by De Masi and Presutti (1991),
the monograph of Spohn (1991),  review papers by
Ferrari (1994, 1996), and the soon-to-appear
 monograph of Kipnis and 
Landim. 

Here are some references that are closely related, but
not directly on the topic of the paper: Hydrodynamic limits
for asymmetric processes with inhomogeneous but not random
rates have been proved by Landim (1996), Covert and 
Rezakhanlou (1997), and Bahadoran (1998). Koukkous (1996)
and Gielis et al.\ (1998) have studied the symmetric ZRP
with random rates. [In the symmetric case the jump
probabilities satisfy $p(x)=p(-x)$.]

\subsubhead Notational remarks
\endsubsubhead
$\mmZ_+=\{0,1,2,3,\ldots\}$. 
 $\ind_A$ and $\ind\{A\}$ denote the indicator random variable
of the event $A$. $\delta_y$ is a delta function or a point mass at $y$,
depending on the context: 
$\delta_y(x)=\ind\{x=y\}$ for points $x$,
and $\delta_x(A)=\ind_A(x)$ for sets $A$.

\head 2. The disordered asymmetric zero-range process \endhead

First we describe  a disordered ZRP
 on $\mmZ^d$ with bounded, monotone jump
 rates. Let $\{p(x):x\in\mmZ^d\}$ be a finite-range probability
distribution, in other words $p(x)\ge 0$, 
$\sum p(x)=1$, and for some fixed finite set $\NN\subseteq\mmZ^d$,
 $p(x)=0$ for $x\notin\NN$. The rate of jumping from a site
depends on the number of particles present through
 a function $r:\mmZ_+\to[0,\infty)$, about which we assume that
$$0=r(0)<r(1)\le r(2)\le r(3)\le r(4)\le \cdots
\tag 2.1
$$
and 
$$r(\infty)=\lim_{k\to\infty}r(k)<\infty.
\tag 2.2
$$
The disorder comes in the form of  random deceleration factors 
$\a_x\le 1$ that depend on the  sites $x$. 
 Once $\mmal=(\a_x:x\in\mmZ^d)$ is picked,
 the dynamics operates as follows: 
 If there
are $\eta(x)\ge 1$ particles at site $x$, then at exponential
 rate $\a_xr(\eta(x))$ a single
particle jumps away from site $x$. The new location
 of this particle is $y$ with probability
$p(y-x)$. This happens at each site $x$ independently of what
happens at other sites. 

 For fixed rates $\mmal$, the generator of the process is  
$$L^{\mmal}f(\eta)=\sum_{x,y\in\mmZ^d} p(y-x)\a_x
r(\eta(x)) [f(\eta^{x,y})-f(\eta)]. 
\tag 2.3
$$
Here $\eta=(\eta(x):x\in{\mmZ^d})$ is an element of the state space
$\SS=\mmZ_+^{\mmZ^d}$
of the process, and $\eta^{x,y}$ is the configuration that results from the 
jump of a single particle from site $x$ to site $y$:
$\eta^{x,y}=\eta+\delta_y-\delta_x$.  Section 5 contains  a
construction of this process, based on a 
percolation argument of Harris (1972). Due to assumption
(2.2) the process can be started from any configuration
$\eta\in\SS$. We denote
the process by 
 $\eta(t)=(\eta(x,t):x\in{\mmZ^d})$, where $t\ge 0$ is the
time variable. 

The standing assumption is that
 $\mmal$ is an ergodic $[c,1]$-valued process for some 
constant $c\in(0,1]$. Let $Q$ denote the distribution of the process
$\mmal$ on the space $\cAA=[c,1]^{\mmZ^d}$.  
Fix $c$ to be the left endpoint of the marginal distribution of $\a_0$, 
so that
 $c$ is the largest number such that the process $\mmal$
is  $[c,1]$-valued. 

What makes the disordered ZRP tractable are 
 invariant distributions that can be explicitly
described. This description uses the same ideas as
 the process without disorder. 
If $\a_x\equiv 1$ (no disorder), among  the extremal
 invariant distributions 
 are the product 
measures $\mu_\psi$ on $\SS$, indexed by a parameter 
$\psi\in[0,r(\infty))$, with  marginals 
$$\mu_\psi\bigl(\eta(x)=k\bigr)=Z(\psi)^{-1}
\;\frac{\psi^k}{r(1)\cdots r(k)}\,,\quad
x\in\mmZ^d\,,\, k\in\mmZ_+\,.
\tag 2.4
$$
 [See Andjel (1982).]
For $k=0$ the product in the denominator is interpreted as 1. 
$Z(\psi)$ is the normalization factor, defined by
$$Z(\psi)=\sum_{k=0}^\infty \frac{\psi^k}{r(1)\cdots r(k)}\,.
\tag 2.5
$$
For the disordered process, fix a choice $\mmal$ for the rates. 
For real numbers $\varphi\in[0,r(\infty)c)$, let
$\nu^\mmal_\varphi$ denote the product probability measure on 
$\SS$ whose marginals vary from site to site, as given by
$$\nu^\mmal_\varphi\bigl(\eta(x)=k\bigr)=Z(\varphi/\a_x)^{-1}
\;\frac{(\varphi/\a_x)^k}{r(1)\cdots r(k)}\,,\quad
x\in\mmZ^d\,,\, k\in\mmZ_+\,.
\tag 2.6
$$

\proclaim{Proposition 1} For each choice of rates $\mmal\in\cAA$
and each value of  $\varphi\in[0,r(\infty)c)$,
the probability distribution $\nu^\mmal_\varphi$ is invariant
for the process with generator $L^\mmal$. 
\endproclaim

The {\it phase transition} of the disordered ZRP is the 
following
situation:  If $Q$ is such that very slow sites are
sufficiently rare, then the family of invariant distributions 
$\{\nu^\mmal_\varphi:\varphi\in[0,r(\infty)c)\}$ does not cover 
the entire range of densities
$0\le \rho<\infty$. Instead, there is a critical density
$\rho^*<\infty$ such that the equilibria $\nu^\mmal_\varphi$ only
exist for densities $\rho\in[0,\rho^*]$. 
To see this, set first
$$M(\psi)=\frac1{Z(\psi)}
\sum_{k=0}^\infty \frac{k\psi^k}{r(1)\cdots r(k)}\,,\quad
\psi\in[0,r(\infty))\,.
\tag 2.7
$$
$M(\psi)$ is the density under $\mu_\psi$.
It is a strictly increasing function from $[0,r(\infty))$ onto
$[0,\infty)$ and has an inverse function $M^{-1}$ which we 
need to refer to below. 
For the disordered model the density $\rho$ as a function of the 
parameter $\varphi$ is defined by averaging over the
random rates: 
$$\rho(\varphi)=
E^Q\biggl[ \frac1{Z(\varphi/\a_0)}
\sum_{k=0}^\infty \frac{k(\varphi/\a_0)^k}{r(1)\cdots r(k)}\biggr]
=E^Q\bigl[ M\bigl(\varphi/\a_0\bigr)\bigr].
\tag 2.8
$$
Here $E^Q$ denotes expectation over the distribution $Q$ of $\mmal$, 
and the random variable inside the expectation is $\a_0$.
For a fixed equilibrium $\nu^\mmal_\varphi$
 the density $\rho(\varphi)$ can be realized as a spatial average
$$\rho(\varphi)=\lim_{\Lambda \nearrow \mmZ^d}
\frac1{|\Lambda|}\sum_{x\in\Lambda}\eta(x)\qquad 
\text{$\nu^\mmal_\varphi$-a.s., for $Q$-a.e.\ $\mmal$.}
$$

By letting $\varphi$ increase to its upper bound
$r(\infty)c$,  (2.8) shows  that the maximal density is 
$$\rho^*=E^Q\bigl[ M\bigl(r(\infty)c/\a_0\bigr)\bigr].
\tag 2.9
$$
This quantity may or may not be infinite, depending on
the distribution $Q$. 
From (2.7) $M(r(\infty))=\lim_{\psi\nearrow r(\infty)}M(\psi)=\infty$,
so in particular if $Q(\a_0=c)>0$, then $\rho^*=\infty$. 
The interesting case with phase transition is the one where $Q(d\a_0)$ has 
a sufficiently thin tail as $\a_0\searrow c$, to make the integral
in (2.9) finite. 

The function $\rho:[0,r(\infty)c)\to[0,\rho^*)$ is
strictly increasing. Let $f:[0,\rho^*)\to[0,r(\infty)c)$ 
denote its inverse function. In other words, for $\rho\in[0,\rho^*)$,
$f(\rho)$ is implicitly defined by
$$\rho=E^Q\bigl[ M\bigl(f(\rho)/\a_0\bigr)\bigr].
\tag 2.10
$$

Now we state a hydrodynamic limit for the disordered ZRP, 
due to Benjamini, Ferrari, and Landim (1996).
For each choice of rates 
$\mmal$, there is  a sequence of
zero-range  processes indexed by $n$, generated by $L^\mmal$. 
$P^\mmal_n$ denotes the probability measure on the probability
space of the $n$th process
  $\{\eta_n(x,t):x\in\mmZ^d, t\ge 0\}$,
$n=1,2,3,\ldots\,$.
 The theorem 
is a weak law of large numbers for the empirical measure
defined by 
$$\pi_n(t)=n^{-d}\sum_{x\in\mmZ^d}\eta_n(x,t)\delta_{x/n},
\tag 2.11
$$
where $\delta_x$ is a unit mass at the point $x\in\mmR^d$.
The assumptions are the following: 

\hbox{}

(A.1) The transition probability $p(x)$ satisfies this
irreducibility condition: for each $x$, $y\in\mmZ^d$ there exists
a finite sequence 
$x=x_0, \dots, x_k=y$ such that $p(x_{i+1} - x_i) 
+ p(x_i - x_{i+1})>0$ for all $i$. 

\hbox{}

(A.2) There exists a bounded continuous function $u_0$ on $\mmR^d$
such that $\|u_0\|_\infty\le 
M(r(\infty)\theta )$ for some $\theta<c$, and 
 for each $\mmal$, 
 the initial distribution of the process $\eta_n$ is
given by 
$$P^\mmal_n\bigl(\eta_n(x,0)=k\bigr)=\mu_{M^{-1}(u_0(x/n))}
\bigl(\eta(x)=k\bigr).$$
[Recall the definitions of $\mu_\psi$ and $M(\psi)$
 from (2.4) and (2.7).] 

\hbox{}

(A.3) The marginal distribution of $\a_0$ is supported by a finite
set: For some $c=c_1<c_2<\cdots<c_m\le 1$, 
$Q\bigl( \a_0\in\{c_1,\ldots,c_m\}\bigr)=1.$

\hbox{}

Assumption (A.2) ensures that for some fixed 
$\varphi\in[0,r(\infty)c)$ and  all $\mmal$,
 all  the initial distributions are
stochastically dominated by $\nu^\mmal_\varphi$.
This is true because, on the $\eta(x)$-marginal 
 $\nu^\mmal_\varphi=\mu_{\varphi/\a_x}$, and this
dominates $\mu_{r(\infty)\theta }$ as long as 
$\varphi/\a_x\ge r(\infty)\theta $, which in turn 
is true for all $\a_x\in[c,1]$ if $\varphi\ge  r(\infty)\theta$. 

Let $\gamma\in\mmR^d$ be the mean drift under $p(x)$:
$\gamma=\sum_{x\in\mmZ^d} xp(x).$
Let $u(x,t)$ on $\mmR^d\times[0,\infty)$
 be the unique entropy solution of the scalar
conservation law
$$\partial u/\partial t+\gamma \cdot\nabla_x [f(u)]=0,\qquad u(x,0)=u_0(x).
\tag 2.12
$$
Let $C_0(\mmR^d)$ denote the space of  compactly
supported continuous functions on $\mmR^d$. 
This theorem was proved by Benjamini et al. (1996): 

\proclaim{Theorem 1} Under assumptions 
{\rm (A.1)--(A.3)}, the following holds 
for $Q$-a.e.\ $\mmal$: For each $t>0$, 
 $\phi\in C_0(\mmR^d)$, and $\e>0$: 
$$\lim_{n\to\infty} 
P_n^\mmal\biggl(\,\biggl|\pi_n(nt,\phi)-\int_{\mmR^d}\phi(x)u(x,t)dx\biggr|
\ge\e\biggr)=0.
\tag 2.13
$$
\endproclaim

The integral against $\pi_n(t)$ is defined by
$\pi_n(t,\phi)=n^{-d}\sum_{x}\eta_n(x,t)\phi(x/n)$. 
The statement is that, in the topology of Radon measures on
$\mmR^d$, $\pi_n(nt)$ converges to $u(x,t)dx$
in probability as $n\to\infty$. 

The shortcoming of this result is that it does not indicate 
what happens on the hydrodynamic scale if the process starts
at density above critical, that is, $u_0(x)>\rho^*$ for some
or all $x$. In fact, assumption (A.3) makes $\rho^*=\infty$,
so there can be no phase transition under these hypotheses. 

Next we state a theorem that covers the hydrodynamics also 
in the high-density regime $\rho>\rho^*$ and admits more
general initial distributions for the process. 
However, we pay a serious price for this strengthening: 
The theorem is valid only for the most basic type
of ZRP with rate function $r(k)=\ind\{k\ge 1\}$.
Furthermore, we are restricted to totally asymmetric jumps
in one dimension:  $d=1$ and $p(1)=1$, so jumps happen
only to the right on $\mmZ$. Finally, we assume that the process
of rates $(\a_x: x\in\mmZ)$ is i.i.d. 

In this case the measures $\nu^\mmal_\varphi$ are products of
geometric distributions:
$$\nu^\mmal_\varphi\bigl(\eta(x)=k\bigr)=(1-\varphi/\a_x)
(\varphi/\a_x)^k\,,\quad
x\in\mmZ\,,\, k\in\mmZ_+\,.
\tag 2.14
$$
Now $M(\psi)=\psi/(1-\psi)$, so the definition of the critical
density becomes
$$\rho^*=c\int_{[c,1]}(\a_0-c)^{-1}Q(d\a_0). 
$$
From this formula it is plainly obvious how the tail of 
 $Q(d\a_0)$ at $\a_0=c+$ determines whether $\rho^*<\infty$
or not, that is, whether phase transition happens or not. 

For $\rho\in[0,\rho^*)$ the flux function $f(\rho)$ is defined
implicitly by the equation 
$$\rho=f(\rho)\int_{[c,1]}\bigl(\a_0-f(\rho)\bigr)^{-1}Q(d\a_0). 
\tag 2.15
$$
If $\rho^*<\infty$, set 
$$\text{$f(\rho)=c$ for $\rho\ge \rho^*$.}
\tag 2.16
$$
This makes $f$ a nondecreasing and concave function on $[0,\infty)$. 

Again we assume we have a sequence of processes $\eta_n(t)$
and corresponding probability measures $P^\mmal_n$. The 
assumption on initial distributions is this: 

\hbox{}

(A.4) Suppose $u_0(x)$ is a nonnegative locally integrable
function on $\mmR$. 
Assume that this holds for $Q$-a.e.\ $\mmal$: 
For all $\phi\in C_0(\mmR)$
and $\e>0$,
$$\lim_{n\to\infty} 
P_n^\mmal\biggl(\,\biggl|\pi_n(0,\phi)-\int_{\mmR}\phi(x)u_0(x)dx\biggr|
\ge\e\biggr)=0.
\tag 2.17
$$

\hbox{}

Previously this assumption was a consequence of
assumption (A.1) so it was not stated explicitly. 
Let $u(x,t)$ on $\mmR\times[0,\infty)$
 be the unique entropy
solution of 
$$\partial u/\partial t+\partial f(u)/\partial x=0,\qquad u(x,0)=u_0(x)
\tag 2.18
$$
where $f$ is the function defined by (2.15)--(2.16). Then we have a theorem 
 due to Sepp\"al\"ainen and Krug (1998):
 
\proclaim{Theorem 2} Assume that $d=1$, $p(1)=1$, 
{\rm $r(k)=\ind\{k\ge 1\}$}, and that $Q$ is 
an  i.i.d.\ distribution for $\mmal=(\a_x)$.
 Then under assumption {\rm (A.4)} 
 the following holds 
for $Q$-a.e.\ $\mmal$: For each $t>0$, 
 $\phi\in C_0(\mmR)$, and $\e>0$: 
$$\lim_{n\to\infty} 
P_n^\mmal\biggl(\,\biggl|\pi_n(nt,\phi)-\int_{\mmR}\phi(x)u(x,t)dx\biggr|
\ge\e\biggr)=0.
\tag 2.19
$$
\endproclaim

In the phase transition case this  theorem has an
 interesting consequence: Suppose the
initial profile satifies $u_0(x)\ge\rho^*$ everywhere on $\mmR$.
Then since $f$ is constant for this range of densities, it follows
that $u(x,t)=u_0(x)$ for all $t>0$. In other
words, the profile does not change 
on the hydrodynamic scale. 

\head 3. The disordered asymmetric exclusion process \endhead

The disordered SEP is less well understood than the
ZRP. The reason is that no invariant distributions have been found. 
Presently we can prove the existence of a hydrodynamic limit
for the totally asymmetric SEP, in dimension one.
 But the flux function 
$f$ of the macroscopic equation (2.18) remains unknown. 
The theorem covers a more general totally asymmetric SEP, namely
the so-called $K$-exclusion, where each site admits $K$ particles
instead of just one. The state space is 
$\SS=\{0,\ldots,K\}^\mmZ$, and we write again
$\eta=(\eta(x):x\in\mmZ)\in\SS$ for the
particle configurations. When the rates $\mmal$ have been chosen, the
generator is
$$L^{\mmal}f(\eta)=\sum_{x\in\mmZ}\a_x
\ind_{\{\eta(x)\ge 1,\eta(x+1)\le K-1\}}
[f(\eta^{x,x+1})-f(\eta)]. 
\tag 3.1
$$
In other words, a jump occurs from site $x$ to $x+1$ at rate 
$\a_x$, provided
site $x$ is not empty and site $x+1$ has less than $K$ particles. 
Consider $K$ fixed but arbitrary. As before, assume we have
a sequence of totally asymmetric $K$-exclusion
processes $\eta_n(t)$, with probability measures $P^\mmal_n$ 
when the rates $\mmal$ are fixed. As for Theorem 2, we only assume
that the initial distributions of the processes have a well-defined
macroscopic profile:

\hbox{}

(A.5) Suppose $u_0$ is a bounded measurable function on $\mmR$ 
such that $0\le u_0(x)\le K$. Assume (2.17) holds for all $\phi\in C_0(\mmR)$
and $\e>0$,  for $Q$-a.e.\ $\mmal$. 

\hbox{}

The theorem about the existence of the hydrodynamic limit is
proved in Sepp\"a\-l\"ai\-nen (1998): 

\proclaim{Theorem 3} Fix a positive integer $K$ and an 
i.i.d.\ distribution $Q$ for $\mmal=(\a_x)$. Then there exists
a concave function $f_K$ on $[0,K]$ that depends on $Q$ such that,
if $u(x,t)$ is the unique entropy solution of {\rm (2.18)} with flux
$f=f_K$, then under assumption {\rm (A.5)},
the limit in {\rm (2.19)} holds for $Q$-a.e.\ $\mmal$, for 
 each $t>0$, 
 $\phi\in C_0(\mmR)$, and $\e>0$. 
\endproclaim

The entropy solution to (2.18) with $f=f_K$  
and $u_0$ bounded measurable can be constructed
 without explicit reference to the p.d.e:
Pick a function $U_0$ such that $U_0'=u_0$ a.e. on $\mmR$. 
Let $f^*_K$ be the concave conjugate of $f_K$:
$$f^*_K(x)=\inf_{0\le \rho\le K}\{x\rho-f_K(\rho)\}.$$
Set
$$U(x,t)=\sup_{y\in\mmR}\bigl\{ U_0(y)+tf^*_K\bigl((x-y)/t\bigr)
\bigr\}.
\tag 3.2
$$
Finally, let $u(x,t)=(\partial/\partial x)U(x,t)$, a derivative 
that is defined at least a.e. Lax (1957) discusses this 
construction of the entropy solution of a scalar conservation law.
See also Section 3.4 in Evans (1998).  
The conjugate $f^*_K$  has a probabilistic meaning in this
context: It is the macroscopic shape of a growth model 
associated to the $K$-exclusion process [Sepp\"al\"ainen (1998)].

\head 4. Open problems  \endhead

\subsubhead Problem 1 \endsubsubhead
Extend Theorem 1 to describe hydrodynamic behavior 
also  at 
high density $\rho>\rho^*$. Or, equivalently, extend
Theorem 2 to more general ZRP's. The proof of Theorem 2 
 in  Sepp\"al\"ainen and Krug (1998) depends 
on special properties of the totally asymmetric 
ZRP with rate function $r(k)=\ind\{k\ge 1\}$. The 
proof of Benjamini et al.\ for Theorem 1 follows
a strategy invented by Rezakhanlou (1991), and may  be 
a better candidate for generalization. 

\subsubhead Problem 2 \endsubsubhead In the phase
transition case, the hydrodynamic 
theorem reveals only trivial behavior at high density $\rho>\rho^*$
(the profile does not move, see Theorem 2). 
It is expected 
that on a finer space scale one can see inhomogeneities develop,
specifically, that the particles concentrate
on exceptionally slow sites. No rigorous results exist
to decribe these phenomena. Krug (1998)
and  Sepp\"al\"ainen and Krug (1998)
discuss this in terms of the exclusion version of the 
totally asymmetric ZRP. 

\subsubhead Problem 3 \endsubsubhead
An interesting open problem for the disordered ZRP
 concerns the weak convergence
of the process: Fix the rates $\mmal$, and take the initial
distribution of $\eta(0)$ spatially ergodic with density $\rho$. 
Does the distribution of  $\eta(t)$ converge weakly, as 
$t\to\infty$, to one of the equilibria $\nu^\mmal_\varphi$?
If $\rho>\rho^*$  in the phase transition case,
does  $\eta(t)$ converge weakly to the equilibrium with
density $\rho^*$? 
For the standard  ZRP results of this type were
proved by Andjel et al.\ (1985). 

\subsubhead Problem 4 \endsubsubhead
For the disordered TASEP,  any new rigorous results
would be welcome. For example, is there a phase transition
similar to the one for disordered ZRP? Does the flux
function $f_K(\rho)$ have a flat segment on an interval
around $\rho=K/2$?

\head 5. Comments on the proofs  \endhead

\subhead 5.1 Construction and invariant distributions \endsubhead
To construct the disordered ZRP with generator (2.3),
start by giving each site $x\in\mmZ^d$ an independent 
rate $r(\infty)$ Poisson point process $\TT^x$ on the time
line $(0,\infty)$. [Here we make use of 
the assumption  $r(\infty)<\infty$.]
 To the $i$th epoch of $\TT^x$ attach two random objects: 
A random threshold $U^x_i$  uniformly distributed
on $[0,r(\infty)]$ and independent of everything else,
and a destination site $y^x_i$ chosen with probability $p(y^x_i-x)$, again 
 independently  of everything else. 
 Fix the rates 
$\mmal=(\a_x:x\in\mmZ^d)$ and an initial configuration
$\eta=(\eta(x):x\in\mmZ^d)$. Informally speaking,
from these ingredients the construction 
of the dynamics goes as follows: Suppose $\tau$ is an
epoch of $\TT^x$ with threshold $U$ and destination site $y$, and the 
dynamics $\eta(t)$ has been constructed for times $0\le t<\tau$. 
If $U\ge \a_xr(\eta(\tau-))$, do nothing. 
If $U< \a_xr(\eta(\tau-))$,  move one particle from 
site $x$ to site $y$; in other words, set 
$$\aligned
\eta(x,\tau)&=\eta(x,\tau-)-1,\\
\eta(y,\tau)&=\eta(y,\tau-)+1,\\
\text{and}\qquad \eta(z,\tau)&=\eta(z,\tau-)\quad\text{ for $z\ne x,y$.}
\endaligned
$$
Subsequently site $x$ lies dormant until the next epoch
of $\TT^x$. However, 
 the lattice $\mmZ^d$ is infinite, so this induction never
even begins because there is no first epoch among 
the point processes $\{\TT^x:x\in\mmZ^d\}$. 

To make the 
construction rigorous, we show that there exists a fixed
time $t_0>0$ such that, starting with an arbitrary $\eta$
in the state space $\SS$, the evolution $\eta(t)$ can be 
computed for $t\in[0,t_0]$. Since $t_0$ is independent
of $\eta$, the construction can be repeated, starting with
$\eta(t_0)$, and extended to $[0,2t_0]$. And so on, to
arbitrarily large times. 

Given a fixed number $t_0>0$
and the Poisson processes $\{\TT^x\}$, construct the following
random graph, with vertex set $\mmZ^d$: Connect $x$ and $y$
with an edge if $x-y$ or $y-x$ lies in $\NN$, and either $\TT^x$ or
$\TT^y$ has an epoch in $[0,t_0]$.  Recall that
$p(z)=0$ for $z$ outside $\NN$. 

\proclaim{Lemma 5.1} For small enough $t_0>0$, the random
graph thus constructed has no infinite connected components,
for almost every realization of $\{\TT^x\}$.
\endproclaim

Before proving the lemma, let us observe how it solves
the construction problem: All the sites $y$ that could
possibly influence the evolution at site $x$ up to time 
$t_0$ are connected to $x$ in the random graph. Since 
$x$ lies in a finite connected component $\CC$, the 
point process $\cup_{x\in\CC}\TT^x$ has only finitely
many epochs in $[0,t_0]$. Consequently the evolution
$\eta(z,t)$ can be computed for $z\in\CC$ and $t\in[0,t_0]$
from the informal rule spelled out above, by considering 
the finitely many epochs in their temporal order. This
procedure is repeated for all connected components. 

\demo{Proof of Lemma 5.1} By translation invariance, it 
suffices to show that the origin is almost surely
connected to only finitely many vertices. Let
$\NN^*=\NN\cup(-\NN)$ and 
$K=|\NN^*|$, the number of vertices $x$ such that either
$x$ or $-x$ lies in $\NN$. Call $x_0,x_1,\ldots,x_n$ a 
self-avoiding path with $n$ edges
 in the random graph
 if $x_i\ne x_j$ for $i\ne j$ and 
there is an edge between $x_i$ and $x_{i+1}$ for each $i$.
In particular, then $x_{i+1}-x_i\in\NN^*$ for each $i$. 

Let $|\cdot|$ be any norm on $\mmR^d$,
and  $R=\max\{|x|:x\in\NN^*\}$. 
If the origin is connected to a vertex $x$ with $|x|>L$,
there is a self-avoiding path with at least $L/R$ 
edges starting at the origin. The probability
that a  self-avoiding path with $2n-1$ edges
starts at the origin is at most
$$K^{2n-1}\bigl(1-e^{-2r(\infty)t_0}\bigr)^n\,.
$$
To see this, note first that the factor $K^{2n-1}$ is an upper
bound on the number of such paths. If 
$0=x_0,x_1,\ldots,x_{2n-1}$ is such  a path, the $n$
edges $(x_0,x_1)$,  $(x_2,x_3)$, $\ldots$,  $(x_{2n-2},x_{2n-1})$ 
are present independently of each other, and each with
probability $1-e^{-2r(\infty)t_0}$ [at least one
of $\TT^{x_{2i}}$ and  $\TT^{x_{2i+1}}$ must have an epoch in
$[0,t_0]$, and each $\TT^{x_{i}}$ has rate $r(\infty)$]. 
Pick $t_0$ small enough so that 
$K^2\bigl(1-e^{-2r(\infty)t_0}\bigr)<1$. Then by Borel-Cantelli,
self-avoiding paths from the origin have a finite
upper bound on their length, almost surely. 
\qed
\enddemo

This approach to the construction of a particle system
is due to Harris (1972). Our presentation followed 
Durrett (1995). 

Let $(\Omega,\FF,\bP)$ denote the probability space 
whose sample point $\omega$ represents a realization
of the Poisson processes $\{\TT^x\}$,  the
random thresholds $\{U^x_i\}$, and the destination 
 sites $\{y^x_i\}$.
We constructed  the random path 
$\eta(\cdot)=\{\eta(x,t):x\in\mmZ^d, t\ge 0\}$ 
 as a function of
the initial state $\eta$, the rates $\mmal$,
and a sample point $\omega$.  Since the Poisson
processes are Markovian, and the random choices of
threshold $U$ and destination state $y$ are
 independent of everything, the process $\eta(\cdot)$
is a time-homogeneous Markov process. 
Let  $\DD=\DD([0,\infty),\SS)$ denote the space of 
right-continuous $\SS$-valued functions on $[0,\infty)$ that
have left limits at each point $t\in[0,\infty)$. 
 For fixed $(\mmal, \eta)$ the path 
 $\eta(\cdot)$ is a $\DD$-valued  random
variable, and it 
 has a probability  distribution $P^{\mmal,\eta}$
on $\DD$ induced from the measure $\bP$ on $\Omega$. 
Write  $E^{\mmal,\eta}$ for the expectation under 
 $P^{\mmal,\eta}$.

Next we prove 
Proposition 1 about the invariant distributions. 
For this we 
restrict the process to a cube 
$\Lambda_k=\{-k,\ldots,k\}^d\subseteq\mmZ^d$. 
The jump probabilities  $p(y-x)$ are then interpreted
with  periodic 
boundary conditions, and become 
$$p_k(x,y)=\sum\{p(z-x): \text{$z=y+(2k+1)w$ for some $w\in\mmZ^d$}\}
\tag 5.1
$$ 
for $x,y\in\Lambda_k$.
The finite-volume generator is 
$$
L^{\mmal}_kf(\eta)=\sum_{x,y\in\Lambda_k} p_k(x,y)\a_x
r(\eta(x)) [f(\eta^{x,y})-f(\eta)]. 
\tag 5.2
$$
$L^{\mmal}_k$ generates
 a Markov jump process with uniformly bounded rates
 on the countable 
state space 
$\SS_k=\mmZ_+^{\Lambda_k}$.  Existence of this process
 follows from standard textbook material. 
Write $E_k^{\mmal,\eta}$ for expectations under distributions
of the process restricted to $\Lambda_k$. 
Notice that if $f$ is a cylinder function  on $\SS$,
$L^{\mmal}_k f=L^\mmal f$ for all large enough $k$. 

Define dual transition probabilities by
$p_k^*(x,y)=p_k(y,x)$ and $p^*(x)=p(-x)$. Then 
$p^*_k$ is obtained from $p^*$ exactly as $p_k$ from $p$
according to (5.1). Let $L^{\mmal,*}_k$ and $L^{\mmal,*}$ be the 
generators with kernels $p^*_k$ and $p^*$ in place of 
 $p_k$ and $p$.

\proclaim{Lemma 5.2} Let $\mmal\in\cAA$ and $\varphi\in[0,r(\infty)c)$,
and let $\nu^\mmal_\varphi$ be the product probability measure
defined by {\rm (2.14)}. Then for all  bounded measurable  
functions $f,g$ on $\SS_k$,
$$\nu^\mmal_\varphi\bigl(gL^\mmal_k f\bigr)=
\nu^\mmal_\varphi\bigl(fL^{\mmal,*}_k g\bigr).
$$
\endproclaim

\demo{Proof}
Start by  writing 
$$\aligned
\nu^\mmal_\varphi\bigl(gL^\mmal_k f\bigr)
&=\sum_{x,y\in\Lambda_k} p_k(x,y)\a_x
\nu^\mmal_\varphi\bigl[r(\eta(x))g(\eta) f(\eta^{x,y})\bigr]\\
&\qquad\qquad -\sum_{x,y\in\Lambda_k} p_k(x,y)\a_x
\nu^\mmal_\varphi\bigl[r(\eta(x))g(\eta) f(\eta)\bigr]\\
&\equiv\ A_1 -A_2,
\endaligned
$$ 
where the last equality defines the quantities  $A_1$ and $A_2$.

Continue first with $A_1$: For any $x,y$ let 
$\nutil$ denote the marginal
 distribution of $\etatil=(\eta(z):z\ne x,y)$.
Abbreviate $R(m)=r(1)\cdots r(m)$. Then for a fixed $x$,
the sum over $y$ in $A_1$ can be expressed as
$$\aligned
&\sum_{y\in\Lambda_k} p_k(x,y)\a_x
\nu^\mmal_\varphi\bigl[r(\eta(x))g(\eta) f(\eta^{x,y})\bigr]\\
=\;&\sum_{y\in\Lambda_k} p_k(x,y)\a_x
\int\nutil(d\etatil)\sum\Sb m_x\ge 1\\ m_y\ge 0\endSb 
Z(\varphi/\a_x)^{-1} Z(\varphi/\a_y)^{-1} R(m_x)^{-1}
  R(m_y)^{-1}\\ 
&\qquad\qquad(\varphi/\a_x)^{m_x} (\varphi/\a_y)^{m_y} 
r(m_x)g\bigl(\etatil,m_x,m_y\bigr) 
f\bigl(\etatil,m_x-1,m_y+1\bigr)\\
=\;&\sum_{y\in\Lambda_k} p_k(x,y)\a_y
\int\nutil(d\etatil)\sum\Sb n_x\ge 0\\ n_y\ge 1\endSb 
Z(\varphi/\a_x)^{-1} Z(\varphi/\a_y)^{-1} R(n_x)^{-1}
  R(n_y)^{-1}\\ 
&\qquad\qquad (\varphi/\a_x)^{n_x} (\varphi/\a_y)^{n_y}
r(n_y)g\bigl(\etatil,n_x+1,n_y-1\bigr) 
f\bigl(\etatil,n_x,n_y\bigr)\\
=\;&\sum_{y\in\Lambda_k} p_k(x,y)\a_y
\nu^\mmal_\varphi\bigl[r(\eta(y))g(\eta^{y,x}) f(\eta)\bigr].
\endaligned
$$
In the above calculation we first used the definition (2.14)
 of $\nu^\mmal_\varphi$. The expectation 
is taken  separately over  the marginal distributions
of $\etatil$, $\eta(x)$, and $\eta(y)$. In the second  
expression,  $m_x$ and $m_y$ are summation 
variables that represent integration over
the distributions of   $\eta(x)$ and $\eta(y)$.
Because $r(m_x)=0$ for $m_x=0$, we sum over $m_x\ge 1$. 
In the second  equality we do a change of variable in
the sum:  $n_x=m_x-1$ and $n_y=m_y+1$.
The last equality is just  the definition (2.14)
 of $\nu^\mmal_\varphi$ again. 

Now sum over $x\in\Lambda_k$ and interchange the summation
indices $x,y$ to get
$$A_1=\sum_{x,y\in\Lambda_k} p_k^*(x,y)\a_x
\nu^\mmal_\varphi\bigl[r(\eta(x))g(\eta^{x,y}) f(\eta)\bigr].
$$
In $A_2$ simply observe that $\sum_y p_k(x,y)=1=\sum_y p_k^*(x,y)$
and write
$$\aligned
A_2&=\sum_{x\in\Lambda_k}\lbrakk
\sum_{y\in\Lambda_k} p_k(x,y)\rbrakk \a_x
\nu^\mmal_\varphi\bigl[r(\eta(x))g(\eta) f(\eta)\bigr]\\
&=\sum_{x,y\in\Lambda_k} p_k^*(x,y) \a_x
\nu^\mmal_\varphi\bigl[r(\eta(x))g(\eta) f(\eta)\bigr].
\endaligned
$$
Combining gives
$$\nu^\mmal_\varphi\bigl(gL^\mmal_k f\bigr)=A_1-A_2=
\nu^\mmal_\varphi\bigl(fL^{\mmal,*}_k g\bigr) 
$$
and the lemma is proved.
\qed
\enddemo

\demo{Proof of Proposition 1} Let $t_0>0$ be the number
chosen in Lemma 5.1. Since $\eta(t)$ is Markovian, invariance
follows if we prove that, for any cylinder function $f$ on $\SS$,
$$\int_\SS E^{\mmal,\eta}\bigl[f(\eta(t))\bigr]\nu^\mmal_\varphi(d\eta)
=\int_\SS f\, d\nu^\mmal_\varphi
\tag 5.3
$$
for all $t\in[0,t_0]$. 

For any $\eta\in\SS$, let $\eta_{\Lambda_k}$ denote its restriction
to $\Lambda_k$. Let $\eta^k(\cdot)$ denote the process in $\SS_k$ 
generated by $L^\mmal_k$. (This of course is {\it not} the
restriction of the process $\eta(\cdot)$ to $\Lambda_k$.)
Taking $g\equiv 1$ in Lemma 5.2 shows that $\nu^\mmal_\varphi$
is invariant for the process $\eta^k(\cdot)$, so we have 
$$\int_{\SS} E^{\mmal,\eta_{\Lambda_k}}_k\bigl[f(\eta^k(t))\bigr]
\nu^\mmal_\varphi(d\eta)
=\int_\SS f\, d\nu^\mmal_\varphi
\tag 5.4
$$
for all $k$ large enough so that $f$ can be regarded as 
a function on $\SS_k$. [Here again we rely on
basic facts of continuous-time Markov chains on countable
state spaces. The whole point is of course that the 
space $\SS$ is not countable, so we need something more to
pass from Lemma 5.2 to  (5.3).] It remains to argue that 
$$
 E^{\mmal,\eta_{\Lambda_k}}_k\bigl[f(\eta^k(t))\bigr]
\to E^{\mmal,\eta}\bigl[f(\eta(t))\bigr]
\tag 5.5
$$
boundedly, as $k\to\infty$, for any fixed $t\in[0,t_0]$ and $\eta\in\SS$. 
For  fixed $(\mmal,\eta)$, we can 
construct the processes $\eta(\cdot)$ and  $\eta^k(\cdot)$
on the same probability space $(\Omega,\FF,\bP)$. We only
need to interpret the destination sites $y^x_i$ ``modulo
$\Lambda_k$'' for the process $\eta^k(\cdot)$.
Then both
$f(\eta^k(t))$ and $f(\eta(t))$ are functions on 
$\Omega$, and both expectations in (5.5) are integrals over 
 $(\Omega,\FF,\bP)$. 
  Since $f$ is bounded,
it suffices to show that 
$$f(\eta^k(t))\to f(\eta(t))\qquad\text{ as $k\to\infty$},
\tag 5.6
$$
almost surely on  $(\Omega,\FF,\bP)$.

 Pick $k_0$ large enough so that 
$f(\eta)$ is completely determined by 
$(\eta(x):x\in\Lambda_{k_0})$. Let $\CC_{k_0}$ be the random 
set of vertices connected to $\Lambda_{k_0}$ in the
random graph discussed in Lemma 5.1. By that lemma, $\CC_{k_0}$
is almost surely finite. But then $f(\eta^k(t))= f(\eta(t))$
as soon as $k$ is large enough so that  $\CC_{k_0}\subseteq\Lambda_k$,
because then all the transitions that affect the value of $f$
are performed identically for  $\eta(\cdot)$ and  $\eta^k(\cdot)$
throughout the time interval $[0,t_0]$. This proves (5.6) and 
thereby (5.5), and then (5.4) turns into (5.3) as  $k\to\infty$. 
\qed
\enddemo

\subhead 5.2 The hydrodynamic limits
\endsubhead 
Theorem 1 is Theorem 3.1 from Benjamini et al.\ (1996). It is
proved by deriving a microscopic version of Kruzkov's entropy
inequalities that characterize the entropy solution of  a
conservation law. This idea for proving hydrodynamic limits
of asymmetric processes is due to Rezakhanlou (1991). 

Theorem 2 follows from Theorem 2 in Sepp\"al\"ainen and
Krug (1998). The discussion  in this paper 
 is formulated for a totally asymmetric 
exclusion process where the random  rates are 
attached to the particles. 
 To obtain results  for the ZRP with site disorder,
interpret the occupation numbers $(\eta_i)$ as the 
gaps (= number of empty sites) between successive
exclusion particles. The proof uses 
 a special construction that works
for the totally asymmetric ZRP with rate function 
$r(k)=\ind\{k\ge 1\}$. Whether the technique can be
somehow generalized to deduce results for other 
choices of $r(k)$ is presently unclear. 

Theorem 3 is proved in Sepp\"al\"ainen (1998).
The approach is similar to the one in  Sepp\"al\"ainen and
Krug (1998). It involves coupling the exclusion process 
 with a countable   collection of realizations
of the same process but with a simple initial 
configuration. This coupling has the property that
a microscopic version of the Lax-Oleinik formula (3.2) holds
almost surely. 

\subsubhead Acknowledgment \endsubsubhead I thank
IME at the University of S\~ao Paulo for fruitful
hospitality, and  Pablo Ferrari,
Joachim Krug, Claudio Landim, Gunter Sch\"utz, and
Sunder Sethuraman for valuable discussions
that contributed to this paper.

\head References \endhead 

\hbox{}

\flushpar  
E. Andjel (1982). Invariant measures
for the zero range process. \ap\ 10 525--547. 
  
\hbox{}

\flushpar 
E. Andjel, C. Cocozza-Thivent, and M. Roussignol (1985).
Quelques compl\'ements  sur le processus de misanthropes et
 le processus ``z\'ero-range''. \aihp\ 21 363-382. 
  
\hbox{}

\flushpar 
C. Bahadoran (1998). Hydrodynamical limit for spatially
heterogeneous simple exclusion process. \ptrf\ 110 287--331.  

\hbox{}

\flushpar  
I. Benjamini, P. Ferrari, and C. Landim (1996). 
Asymmetric conservative processes with random
rates. \spa\ 61 181--204. 

\hbox{}

\flushpar  
P. Covert and F. Rezakhanlou (1997). Hydrodynamic limit 
for particle systems with nonconstant speed parameter.
\jsp\ 88 383--426.

\hbox{}

\flushpar  
A. De Masi and E. Presutti (1991).  
 Mathematical Methods for Hydrodynamic Limits. 
Lecture Notes in Mathematics 1501,  
 Springer-Verlag,   Berlin. 
    
\hbox{}

\flushpar  
R. Durrett (1988).
 Lectures on Particle Systems and Percolation. 
Wadsworth and Brooks/Cole. 
   
\hbox{}

\flushpar  
R. Durrett (1995).
Ten  lectures on particle systems. Lecture Notes
in Mathematics 1608 (Saint-Flour, 1993), 97--201. Springer-Verlag.

\hbox{}

\flushpar  
L. C. Evans (1998). Partial Differential Equations. 
American Mathematical Society. 
 
\hbox{}

\flushpar  
P. A. Ferrari (1994). 
Shocks in one-dimensional processes 
with drift. 
Probability and Phase Transitions,
ed. G. Grimmett,  
 Kluwer Academic Publishers,  35--48.

\hbox{}

\flushpar  
P. A. Ferrari (1996). Limit theorems for
tagged particles. 
Markov Process.\ Related Fields 2  17--40.
  
\hbox{}

\flushpar  
G. Gielis, A. Koukkous, and C. Landim (1998). 
Equilibrium fluctuations for zero range processes 
in random environment. \spa\ 77 187--205.

\hbox{}

\flushpar  
D. Griffeath (1979).  Additive
and Cancellative Interacting Particle Systems. 
 Lecture Notes in Mathematics 724, Springer-Verlag. 

\hbox{}

\flushpar  
T. E. Harris (1972). Nearest-neighbor Markov interaction
processes on multidimensional lattices. Adv.\ Math.\ 9
66--89. 
  
\hbox{}

\flushpar  
A. Koukkous (1996). Hydrodynamic behavior of symmetric 
zero range process with random rates. Preprint. 

\hbox{}

\flushpar  
J. Krug (1998). Platoon formation as a critical phenomenon. 
Traffic and granular flow `97, ed. D. E. Wolf and M.
Schreckenberg, Springer-Verlag. 

\hbox{}

\flushpar  
J. Krug and P. Ferrari (1996).
Phase transitions in driven diffusive systems
with random rates.  J. Phys. A: Math. Gen.
 29 L465--L471. 

\hbox{}

\flushpar  
C. Landim (1996). Hydrodynamical limit
for space inhomogeneous one-dimensional totally asymmetric
zero-range processes. \ap\ 24  599--638. 

\hbox{}

\flushpar  
P. Lax (1957) Hyperbolic systems of
conservation laws II. 
Comm.\ Pure Appl.\ Math.  10 537--566.

\hbox{}

\flushpar  
T. M. Liggett (1985). 
  Interacting Particle Systems.
Springer-Verlag,  New York. 

\hbox{}

\flushpar  
F. Rezakhanlou (1991).
Hydrodynamic limit for attractive particle systems
on $\mmZ^d$. Comm. Math. Phys. 
  140  417--448.

\hbox{}

\flushpar  
T. Sepp\"al\"ainen (1998). Existence of hydrodynamics for the 
totally asymmetric simple $K$-exclusion process. 
To appear in \ap

\hbox{}

\flushpar  
T. Sepp\"al\"ainen and  J. Krug (1998). 
Hydrodynamics and platoon formation for a
 totally asymmetric exclusion
model with particlewise disorder. To appear in \jsp
     
\hbox{}

\flushpar  
H. Spohn (1991). 
 Large Scale Dynamics of Interacting Particles.
Springer-Verlag,   Berlin.

\enddocument